\documentclass[11pt]{article}
\usepackage{color,latexsym,amsmath,amssymb,path,url,enumerate,cases}

\newtheorem{theorem}{Theorem}
\newtheorem{theorem*}{Theorem}
\newtheorem{lemma}[theorem]{Lemma}

\newcommand{\Proof}[1]
        {
        \noindent
        \emph{Proof #1.}~
        }
\newsavebox{\smallProofsym}                     
\savebox{\smallProofsym}                        %
        {
        \begin{picture}(7,7)                    %
        \put(0,0){\framebox(6,6){}}             %
        \put(0,2){\framebox(4,4){}}             %
        \end{picture}                           %
        }                                       %
\newcommand{\smalleop}[1]
        {
        \mbox{} \hfill #1~~\usebox{\smallProofsym}\!\!\!\!\!\!\
        }
\newenvironment{theProof}[1]
        {
        \Proof{#1}}{\smalleop{}
        \medskip

        }
\usepackage{dsfont}

\newcommand{\RR}{\ensuremath{\mathbb R}}

\begin{document}
\pagenumbering{arabic}

\title{Distinct distances between a collinear set \\
and an arbitrary set of points\thanks{%
Work on this paper was partially supported by Grant 2012/229 from
the U.S.-Israel Binational Science Foundation, by Grant 892/13 from
the Israel Science Foundation, by the Israeli Centers for Research
Excellence (I-CORE) program (center no.~4/11), by the Blavatnik Computer
Science Research Fund at Tel Aviv University, and by the Hermann
Minkowski--MINERVA Center for Geometry at Tel Aviv University.}}

\author{
Ariel Bruner\thanks{%
Blavatnik School of Computer Science, Tel Aviv University,
Tel Aviv, Israel; \tt{arielbro@mail.tau.ac.il} } \and
Micha Sharir\thanks{%
Blavatnik School of Computer Science, Tel Aviv University, Tel~Aviv,
Israel; \tt{michas@post.tau.ac.il} } 
}

\maketitle

\begin{abstract}
We consider the number of distinct distances between two finite sets of points in $\RR^k$,
for any constant dimension $k\ge2$, where one set $P_1$ consists of $n$ points on a line $l$,
and the other set $P_2$ consists of $m$ arbitrary points, such that no hyperplane orthogonal to $l$
and no hypercylinder having $l$ as its axis contains more than $O(1)$ points of $P_2$.
The number of distinct distances between $P_1$ and $P_2$ is then
$$
\Omega\left(\min\left\{ n^{2/3}m^{2/3},\; \frac{n^{10/11}m^{4/11}}{\log^{2/11}m},\; n^2,\; m^2\right\}\right) .
$$
Without the assumption on $P_2$, there exist sets $P_1$, $P_2$ as above, with only $O(m+n)$ distinct distances
between them.
\end{abstract}

\section{Introduction}

The problem of bounding the minimum number of distinct distances that can be determined by a set of $n$ points
in the plane was first proposed by Paul Erd\H{o}s in 1946 \cite{erdos}. The problem was stated as follows:
\begin{quote}
\emph{``}Let $[P_n]$ be the class of all planar subsets $P_n$ of $n$ points and denote by $f(n)$
the minimum number of different distances determined by its $n$ points for $P_n$ an element of $[P_n]$.\emph{"}
\end{quote}
Erd\H{o}s then proved the bounds $\sqrt{n-3/4}-1/2\leq f(n)\leq c n/\sqrt{\log n}$, stating before the proof:
\begin{quote}
\emph{``}The following theorem establishes rough bounds for arbitrary $n$. Though I have sought
to improve this result for many years, I have not been able to do so.\emph{"}
\end{quote}

Erd\H{o}s conjectured that the actual value of $f(n)$ is much closer to the upper bound. Several improvements have been made to the lower bound
over the years, culminating in the recent bound $f(n)\in \Omega(\frac{n}{\log n})$ of Guth and Katz in 2010 \cite{GuthKatz}. There still
remains a small gap between Erd\H{o}s's conjectured bound and this best known lower bound. Many interesting variants of the problem are still wide
open (that is, with significant gaps between the best known upper and lower bounds),
such as the variant considering point sets in $\RR^d$ for $d\geq3$ \cite{erdos}, or variants in which the points are further constrained;
e.g., so that no three points are collinear \cite{ErdosCollinear}.

One such family of problems concern distinct distances between two sets of points, $P_1$ and $P_2$, i.e.,
the problem is to obtain a lower bound on the cardinality of the set
$$
\mathcal{D}(P_1,P_2) = \left\{\|p-q\|\mid p\in P_1,q\in P_2 \right\} .
$$
One of the most basic variants in this family is a planar instance, where $P_1$ lies on a line $l_1$ and $P_2$ lies on a line $l_2$,
such that $l_1$ and $l_2$ are neither parallel nor orthogonal. Purdy conjectured that in this case, for $|P_1|=|P_2|=n$,
one always has $\mathcal{D}(P_1,P_2) = \omega(n)$ (e.g., see \cite[Section~5.5]{BrassMoserPach}).
Elekes and R\'{o}nyai \cite{ElekesLajos} have proven the conjecture in 2000, without stating a specific bound.
A bound of $\Omega(n^{5/4})$ was later given by Elekes \cite{Elekes}, and was recently improved by
Sharir, Sheffer, and Solymosi to $\Omega(n^{4/3})$ \cite{SharirShefferSolymosi}.
They have also derived a bound for the unbalanced case, where $P_1$ and $P_2$ have different cardinalities.

More recently, Pach and de Zeeuw \cite{PachZeeuw} generalized the result of \cite{SharirShefferSolymosi},
with the same lower bound, for any
pair of sets $P_1$, $P_2$ that lie on two respective (possibly coinciding) irreducible plane algebraic curves
of constant degree, excluding the cases of parallel lines, orthogonal lines, or concentric circles.
In higher dimensions, Charalambides~\cite{Charalambides} obtained the lower bound $\Omega(n^{5/4})$ for
the number of distinct distances between $n$ points on a constant-degree algebraic curve, in any dimension,
unless the curve is an \emph{algebraic helix} (see~\cite{Charalambides} for the definition). The bound has recently
been improved to $\Omega(n^{4/3})$ by Raz~\cite{Raz}.

In this paper we consider the bipartite distinct distances problem in $\RR^k$, for any fixed
dimension $k\ge 2$, involving two sets $P_1$, $P_2$,
where $P_1$ is contained in a line $l$, and $P_2$ is fairly arbitrary. Our main result is the following theorem.

\begin{theorem}\label{th:master_no_polyk}
Let $l$ be a line in $\RR^k$, for any fixed $k\ge 2$, let $c$ be an absolute constant, and let
$P_1$, $P_2$ be two finite point sets of respective cardinalities $n$, $m$, such that $P_1\subset l$ and $P_2$
is an arbitrary set of points in $\RR^k$, such that no hyperplane orthogonal to $l$, and no circular
hypercylinder whose axis is $l$ contains more than $c$ points of $P_2$. That is, assuming $l$ to be the $x_1$-axis,
\begin{align*}
\forall t\in \RR, & \quad |\{p \in P_2 \mid p_1=t\}|\leq c , \quad \text{and} \\
\forall t\in \RR, & \quad |\{p \in P_2 \mid \rho(p,l)=t\}|\leq c ,
\end{align*}
where $\rho(p,l)$ is the (Euclidean) distance of $p$ from $l$, and $p_1$ is the first coordinate of $p$. Then 
\begin{equation}\notag
|\mathcal{D}(P_1,P_2)| \in
\left\{
\begin{array}{lcc}
\Omega (m^2)                                       & \text{if} &  m \leq n^{1/2} \\
\Omega (n^{2/3}m^{2/3})                            & \text{if} & n^{1/2}<m\leq n^{4/5}/\log^{3/5}n \\
\Omega \left(\frac{n^{10/11}m^{4/11}}{\log^{2/11}m} \right)    & \text{if} & n^{4/5}/\log^{3/5}n <m\leq n^3\\
\Omega (n^2)                                       & \text{if} & n^3<m .\\
\end{array}
\right.
\end{equation}
\end{theorem}
The conditions on $P_2$ are necessary---one can easily construct a set $P_1$ of $n$ points on a line, and a set
$P_2$ of $m$ points lying in any of the forbidden hyperplanes or hypercylinders, such that
$|\mathcal{D}(P_1,P_2)| = O(m+n)$. Such constructions are well known, already for the planar case
(e.g., see~\cite{SharirShefferSolymosi}).

So far, distinct distances in higher dimensions have been studied, within the algebraic framework
pioneered by Guth and Katz~\cite{GuthKatz}, in only a few papers, including those mentioned above.
Although these results do not directly apply to the setup considered in this paper, where only one set is contained in a line,
we nevertheless show that this general approach can be adapted to the scenario considered here too.
Our bound is worse than that of Pach and de Zeeuw~\cite{PachZeeuw}, when $|P_2|$ is large relative to $|P_1|$
(we get the same bound for smaller-size sets $P_2$), but the setup in which it applies is considerably more general,
in its assumptions concerning $P_2$. Our result can also be compared with Raz's bound~\cite{Raz}, which applies to general algebraic curves
in higher dimensions, and is the same as in~\cite{PachZeeuw}, but it only handles the non-bipartite case,
and requires all the points to lie on the curve.


\section{Proof of \protect{Theorem~\ref{th:master_no_polyk}}}

\subsection{Distinct distances in the plane}\label{sc:plane}

We first establish Theorem~\ref{th:master_no_polyk} in the plane. The proof in arbitrary dimensions is essentially
the same, except that the notations are somewhat more tedious. After handling the planar case, we will comment
on the minor modifications that are needed for the general case.

Let $l$ be a line in the plane, which, without loss of generality, we take to be the $x$-axis.
Let $P_1$ be a set of $n$ points on $l$, and let $P_2$ be an arbitrary set of $m$ points, such that
no line parallel or orthogonal to $l$ contains more than $c$ points of $P_2$, for some (arbitrary) absolute constant $c$.
Define $\mathcal{D}(P_1,P_2)$ as in the introduction, and put $x=|\mathcal{D}(P_1,P_2)|$.
To establish the lower bound on $x$ asserted in the theorem, we proceed as follows.

Let $P'_2$ be a maximal subset of $P_2$ such that no two points of $P'_2$ have the same $x$-coordinate
or the same $y$-coordinate, and all the points of $P'_2$ lie on the same side of, say, above the $x$-axis.
By the assumptions on $P_2$, we have $|P'_2|\geq \lfloor \frac{m}{2(2c-1)} \rfloor \in \Theta(m)$.
Indeed, we first assume, without loss of generality, that the upper halfplane contains at least
$\lfloor m/2\rfloor$ points of $P_2$, and ignore the points lying below $l$. We then construct $P'_2$
greedily, picking at each step an arbitrary point from $P_2$ and discarding the at most $2(c-1)$ other points with the
same $x$-coordinate or the same $y$-coordinate. We repeat this step until (the subset in the upper halfplane of) $P_2$ is exhausetd, and get a
subset $P'_2$ with all the desired properties.
Since $|\mathcal{D}(P_1,P_2)|\geq|\mathcal{D}(P_1,P'_2)|$, a lower bound for the latter quantity will
serve as a lower bound for the former as well.
To simplify the notation, we continue to refer to the pruned set $P'_2$ as $P_2$, and to its cardinality as $m$.

We follow a standard approach, already used in several earlier works, such as \cite{GuthKatz, SharirShefferSolymosi},
and consider the set
\begin{equation}\notag
Q=\{(a,p,b,q)\mid a,b\in P_1,p,q\in P_2, \|a-p\|=\|b-q\|,(a,p)\neq(b,q)\}
\end{equation}
of quadruples. Enumerate the elements of $\mathcal{D}(P_1,P_2)$ as $d_1,d_2,\dotsc,d_x$, and define
$$
E_i=\Big\{(a,p)\mid a\in P_1,b\in P_2,\|a-p\|=d_i\Big\} ,
$$
for $i=1,\ldots,x$.
Then, by the Cauchy-Schwarz inequality, it follows that
\begin{equation}\notag
|Q| = 2\sum_{i=1}^{x}\binom{|E_i|}{2} \ge
\sum_{i=1}^{x} \left(|E_i|-1 \right)^2 \ge \frac{1}{x}\left(\sum_{i=1}^{x} (|E_i|-1) \right)^2 =
\frac{\left(mn -x\right)^2}{x}.
\end{equation}
We may assume $x\leq \frac{mn}{2}$, for otherwise we are done. Then
\begin{equation}\label{eq:lowerXb}
x\geq \frac{m^2n^2}{4|Q|}.
\end{equation}

To obtain the asserted lower bound on $x$, we now derive an upper bound for $|Q|$.
For that, we employ a result due to Agarwal et al.~\cite{Agarwal}. In the language of \cite{Agarwal},
a family $C$ of \emph{pseudo-parabolas} is a family of graphs of everywhere defined continuous functions,
so that each pair intersect in at most two points. The collection $C$ admits a
\emph{3-parameter algebraic representation} if it is invariant to translation, its elements are all
semi-algebraic with constant description complexity, and have \emph{three degrees of freedom},
in the sense that each curve in $C$ can be specified by three real parameters, and is thus
identifiable with a point in $\RR^3$, in a suitable manner. A full definition of this property
is given at the beginning of \cite[Section 5]{Agarwal}.

\begin{theorem}{\bf{(Agarwal et al., \cite[Theorem 6.6]{Agarwal})}}\label{th:Agarwal}
Let $C$ be a finite family of distinct pseudo-parabolas that admits a 3-parameter algebraic representation,
and let $P$ be a finite set of distinct points in the plane. Then
\begin{equation}\notag
I(P,C)\in O\left(|P|^{2/3}|C|^{2/3}+|P|^{6/11}|C|^{9/11}\log^{2/11}|C|+|P|+|C| \right).
\end{equation}
\end{theorem}
(The bound in \cite{Agarwal} is slightly weaker; the improvement, manifested in the factor $\log ^{2/11}|C|$,
which replaces a slightly larger factor in \cite{Agarwal}, is due to Marcus and Tardos \cite{MarcusTardos}.)

Define $Q_0=\{(a,p,b,q)\in Q\mid p=q\}$ and $Q_1=Q\setminus Q_0$. For each $a\in P_1$, $p\in P_2$, there can be at most one point $b\neq a$ in
$l$ such that $|ap|$=$|bp|$, since a circle around $p$ intersects $l$ in at most two points. Thus,
\begin{equation}\label{eq:Q0}
|Q_0| \in O(nm).
\end{equation}

We provide an upper bound for $|Q_1|$ by reducing the problem to a curve-point incidence problem and then applying Theorem \ref{th:Agarwal}.
Combining this with \eqref{eq:Q0} yields an upper bound for $|Q|$.

Let $a=(x,0),b=(y,0)\in P_1$, and $p=(p_1,p_2),q=(q_1,q_2)\in P_2$, such that $(a,p,b,q)\in Q_1$. Then
\begin{equation}\label{eq:prehyperbola}
\sqrt{(x-p_1)^2 + p^2_2} = \sqrt{(y-q_1)^2 + q^2_2},
\end{equation}
or
\begin{equation}\label{eq:hyperbola}
x^2 -y^2 - 2p_1x + 2q_1y + p^2_1 - q^2_1 + p^2_2 - q^2_2 = 0.
\end{equation}

For fixed $p,q$, equation \eqref{eq:hyperbola} defines a hyperbola $h_{pq}$.
Define $\Gamma = \left\{h_{pq}\mid p\ne q\in P_2 \right\}$.
We claim that the hyperbolas of $\Gamma$ are all distinct. Indeed, let $(p,q)\neq(p',q')$
be two distinct pairs of points from $P_2$. Then either $p\neq p'$ or $q\neq q'$, and, from the assumption that
each point in $P_2$ has a distinct $x$-coordinate, either the coefficients of $x$
or the coefficients of $y$ in $h_{pq}$ and $h_{p'q'}$ must differ, from which it easily follows that the
two hyperbolas are indeed distinct.

Define $\Pi=\{(x,y)\mid (x,0),(y,0)\in P_1\}$. A quadruple $(a=(a_x,0),p,b=(b_x,0),q)$ is in $Q_1$
if and only if $(a_x,b_x)\in \Pi$ is incident to $h_{pq} \in \Gamma$, and thus
\begin{equation}\label{eq:incidences}
|Q_1|=I(\Pi,\Gamma) .
\end{equation}

Two distinct hyperbolas may share a common component only if they are both degenerate (i.e., the product of two lines).
Fortunately, we have the following property.

\begin{lemma}\label{lm:reducible}
There are no degenerate hyperbolas in $\Gamma$.
\end{lemma}
\begin{theProof}{}
By \eqref{eq:prehyperbola}, the hyperbola $h_{pq}$ can be written as
\begin{equation}\notag
(x-p_1)^2 - (y-q_1)^2 = q_2^2 - p_2^2,
\end{equation}
which can degenerate into two lines if and only if $q_2^2 = p_2^2$.
Since we assume that the points of $P_2$ have distinct $y$-coordinates, all
of the same sign, we have $q_2^2 \neq p_2^2$ (recall that we only consider here hyperbolas $h_{pq}$ with $p\neq q$), and the lemma follows.
\end{theProof}

All hyperbolas of $\Gamma$ are contained in the family $H$ of hyperbolas whose equations are of the form $(x+\alpha)^2-(y+\beta)^2 + \gamma=0$,
where $\gamma$ is either positive or negative (but not zero).
Note that half of the hyperbolas of $\Gamma$ satisfy $\gamma>0$ and the other half satisfy $\gamma<0$
(one is obtained from the other by interchanging $p$ and $q$).
The treatments of these two subclasses of hyperbolas are essentially identical, and we focus here on those
hyperbolas with $\gamma > 0$, and continue to denote this set as $\Gamma$. Each of these hyperbolas has a
top branch and a bottom branch. We partition each hyperbola in $H$ into its two branches, and denote by
$H_T$ (resp., $H_B$) the set of the top (resp., bottom) branches of the hyperbolas in $H$.
Each branch is a semialgebraic set of degree 2, and the graph of a continuous totally defined function.
Moreover, as is easily seen, each pair of elements of $H_T$ or of $H_B$ intersect in at most two points,
and these families are invariant under translation. Finally, each of the branches in $H_T$ or in $H_B$
can be specified by the three real parameters $\alpha, \beta, \gamma$. Hence, $H_T$ and $H_B$ meet the
conditions in Theorem \ref{th:Agarwal}. Let $\Gamma_T$ (resp., $\Gamma_B$) denote the set of the top
(resp., bottom) branches of the hyperbolas of $\Gamma$. We have $\Gamma_T \subseteq H_T$ and $\Gamma_B \subseteq H_B$.
Combining Theorem \ref{th:Agarwal} with \eqref{eq:incidences}, we have
\begin{equation}\notag
I(\Pi,\Gamma_L),I(\Pi,\Gamma_R)\in O\Big(|\Pi|^{2/3}|\Gamma|^{2/3}+|\Pi|^{6/11}|\Gamma|^{9/11}\log^{2/11}|\Gamma|+|\Pi|+|\Gamma|\Big) ,
\quad\text{and thus}
\end{equation}
\begin{equation}\notag
|Q_1| = I(\Pi,\Gamma_L) + I(\Pi,\Gamma_R) \in O\Big(n^{4/3}m^{4/3}+n^{12/11}m^{18/11}\log^{2/11}m+n^2+m^2\Big) .
\end{equation}
Since this bound subsumes the bound for $|Q_0|$ in \eqref{eq:Q0}, it is also
an upper bound for $|Q|$. Combining this bound with \eqref{eq:lowerXb}, we then get
\begin{equation}\notag
x \in
\left\{
\begin{array}{lcc}
\Omega (m^2)                                       & \mbox{if} &  m \leq n^{1/2} \\
\Omega (n^{2/3}m^{2/3})                            & \mbox{if} & n^{1/2}<m\leq n^{4/5}/\log^{3/5}n \\
\Omega \left(\frac{n^{10/11}m^{4/11}}{\log^{2/11}m} \right)    & \mbox{if} & n^{4/5}/\log^{3/5}n <m\leq n^3\\
\Omega (n^2)                                       & \mbox{if} & n^3<m ,\\
\end{array}
\right.
\end{equation}
thereby completing the proof.

\subsection{Distinct distances in higher dimensions}\label{sc:space}

In this subsection we establish the general version of Theorem \ref{th:master_no_polyk}, in higher dimensions, by modifying the proof
for the planar case in a fairly straightforward manner. Let $k>2$ be some fixed dimension, let $l$ be a line in $\RR^k$,
which, without loss of generality, we take to be the $x_1$-axis. Let $P_1$ be a set of $n$ points on $l$,
and let $P_2$ be an arbitrary set of $m$ points in $\RR^k$, such that no hyperplane orthogonal to $l$, and no circular
hypercylinder whose axis is $l$ contains more than $c$ points of $P_2$, for some fixed constant parameter $c$. That is,
\begin{align}
\label{eq:Ccond1}
\forall t\in \RR, & \quad |\{p \in P_2 \mid p_1=t\}|\leq c , \quad \text{and} \\
\forall t\in \RR, & \quad |\{p \in P_2 \mid \rho(p,l)=t\}|\leq c ,
\label{eq:Ccond2}
\end{align}
Notice that this extends the assumption made in the planar case, that no line parallel or orthogonal
to $l$ contains more than $c$ points of $P_2$.

As in Subsection \ref{sc:plane}, we can assume that $P_2$ satisfies \eqref{eq:Ccond1} and \eqref{eq:Ccond2}
with $c=1$. This is achieved by the same greedy construction, repeatedly choosing an arbitrary point $p$ from $P_2$, and then discarding all the
at most $2(c-1)$ other points that lie in the same hyperplane orthogonal to $l$ or in the same circular hypercylinder with $l$ as an axis.
We get a subset $P'_2$ of $P$ with at least $\lfloor m/(2c-1)\rfloor = \Theta(m)$ points with the desired property. Continue
to denote $P'_2$ as $P_2$ and its size as $m$.

The proof proceeds along the same lines as the previous proof. We consider the set $Q$ of quadruples, partition it into
the subsets $Q_0$ and $Q_1$, where the quadruples $(a,p,b,q)$ in $Q_0$ (resp., in $Q_1$) satisfy $p=q$ (resp., $p\ne q$).
Since a hypersphere intersects $l$ in at most two points, we again have
\begin{equation}\label{eq:Q0k}
|Q_0| \in O(nm).
\end{equation}
For each quadruple $(a,p,b,q)\in Q_1$, we have
\begin{equation}\notag
\sqrt{(a_1-p_1)^2 + \sum_{i=2}^{k}{p^2_i}} = \sqrt{(b_1-q_1)^2 + \sum_{i=2}^{k}{q^2_i}},
\end{equation}
or
\begin{equation}\label{eq:prehyperbolak}
a_1^2 -b_1^2 - 2p_1a_1 + 2q_1b_1 + \left( \sum_{i=1}^{k}{p^2_i} - \sum_{i=1}^{k}{q^2_i} \right)= 0.
\end{equation}
As before, for a given pair $(p,q)$, the set of pairs of $x_1$-coordinates of all $a,b$ that satisfy
\eqref{eq:prehyperbolak} is a hyperbola, defined by the equation

\begin{equation}\label{eq:hyperbolak}
(x-p_1)^2 -(y-q_1)^2 + \left( \sum_{i=2}^{k}{p^2_i} - \sum_{i=2}^{k}{q^2_i} \right)= 0.
\end{equation}

As before, we denote this hyperbola as $h_{pq}$. Arguing as in the planar case (using \eqref{eq:Ccond1} with $c=1$),
we conclude that two distinct pairs of points of $P_2^2$ cannot define the same hyperbola.
As in Section \ref{sc:plane}, a quadruple $(a,p,b,q)$ is in $Q_1$ if and only if the point
$(a_1,b_1)$ is incident to the hyperbola $h_{pq}$, and all these hyperbolas are non-degenerate
(because of \eqref{eq:Ccond2} with $c=1$). Thus,

\begin{equation}\label{eq:incidencesk}
|Q_1|=|\mathcal{I}(\Pi,\Gamma)|.
\end{equation}

$\Gamma$ is contained in the same family of hyperbolas defined in the last section, namely
hyperbolas that can be defined by $(x+\alpha)^2-(y+\beta)^2 + \gamma=0$. Restricting our attention
to hyperbolas that have top and bottom branches (those with $\gamma > 0$), and
treating the sets $\Gamma_T$, $\Gamma_T$ of top branches and bottom branches of these hyperbolas
separately, the same analysis then yields
\begin{equation}\notag
I(\Pi,\Gamma_L),I(\Pi,\Gamma_R)\in O\Big(|\Pi|^{2/3}|\Gamma|^{2/3}+|\Pi|^{6/11}|\Gamma|^{9/11}\log^{2/11}(|\Gamma|)+|\Pi|+|\Gamma|\Big),
\quad\text{and thus}
\end{equation}
\begin{equation}\notag
|Q_1|=I(\Pi,\Gamma_L)+I(\Pi,\Gamma_R)\in O\Big(n^{4/3}m^{4/3}+n^{12/11}m^{18/11}\log^{2/11}m+n^2+m^2\Big) .
\end{equation}
Since this bound subsumes the bound for $|Q_0|$ obtained in \eqref{eq:Q0k}, it is also an upper bound for $|Q|$, and so
\begin{equation}\notag
|Q|\in O\left(n^{4/3}m^{4/3}+n^{12/11}m^{18/11}\log^{2/11}m+n^2+m^2\right).
\end{equation}

Then, by combining this bound with \eqref{eq:lowerXb}, we get
\begin{equation}\notag
x \in
\left\{
\begin{array}{lcc}
\Omega (m^2)                                       & \text{if} &  m \leq n^{1/2} \\
\Omega (n^{2/3}m^{2/3})                            & \text{if} & n^{1/2}<m\leq n^{4/5}/\log^{3/5}n \\
\Omega \left(\frac{n^{10/11}m^{4/11}}{\log^{2/11}m} \right)    & \text{if} & n^{4/5}/\log^{3/5}n <m\leq n^3\\
\Omega (n^2)                                       & \text{if} & n^3<m ,\\
\end{array}
\right.
\end{equation}
completing the proof of the theorem.
$\Box$

\paragraph{Remarks.}
An obvious open problem is to extend the result to the case where $P_1$ lies on some constant-degree algebraic curve,
and $P_2$ is arbitrary. This is probably doable, except that the hyperbolas $h_{pq}$ are now replaced by more general
algebraic curves, which no longer satisfy the pseudo-parabolas properties. We leave it as a topic for further research
to explore such an extension.

In addition, the requirement that every hyperplane or hypercylinder of the above kinds contains at most $O(1)$ points of
$P_2$ is perhaps rather restrictive. What happens if this number is constrained to be at most $s$, for some
non-constant $s\ll |P_2|$? It would be interesting to find a sharp dependence of the number of distinct distances on $s$.

\paragraph{Acknowledgements.}
The authors, especially the first one, would like to thank Adam Sheffer for valuable discussions and for his guidance.


\end{document}